\documentclass[smallextended,numbook,runningheads]{svjour3}     
\smartqed  

\usepackage{graphicx}
\usepackage{amsmath}
\usepackage{epstopdf} 
\usepackage{mathptmx}      

\usepackage{algorithm}
\usepackage{algpseudocode}

\usepackage{booktabs}
\usepackage{url}

\renewcommand{\Re}{\hbox{I\hskip -1.8pt R}}
\DeclareMathOperator*{\argmax}{argmax} 
\DeclareMathOperator*{\argmin}{argmin} 

\usepackage{color}

\journalname{BIT}
\begin{document}

\title{Solving large linear least squares problems with linear equality constraints\thanks{The second author was supported
by  project GAGA-12719S  of the Grant Agency of the Czech Republic.
}}

\titlerunning{Linear least squares problems with linear equality constraints}

\author{Jennifer Scott         \and
        Miroslav T\accent23uma
}

\institute{Jennifer Scott \at
              STFC Rutherford Appleton Laboratory,
Harwell Campus, Didcot, Oxfordshire, OX11 0QX, UK
and School of Mathematical, Physical and Computational Sciences,
University of Reading, Reading RG6 6AQ, UK. \\
              \email{jennifer.scott@stfc.ac.uk}           
           \and
           Miroslav T\accent23uma \at
Department of Numerical Mathematics, Faculty of Mathematics and Physics,
Charles University, Czech Republic.\\
              \email{mirektuma@karlin.mff.cuni.cz}
              }
\date{Received: date / Accepted: date}

\maketitle

\begin{abstract}
We consider the problem of efficiently solving large-scale linear
least squares problems
that have one or more linear constraints that must be satisfied exactly.
Whilst some classical approaches are theoretically well founded, they can face
difficulties when the matrix of constraints contains dense rows or if an algorithmic transformation
used in the solution process results in a modified problem that is much denser than the original one.
To address this, we propose modifications and new ideas,
with an emphasis on requiring that the constraints be satisfied with a small residual.
We examine combining the null-space method  with our recently developed
algorithm for computing a null space basis matrix for a ``wide'' matrix. We further show that
a direct elimination approach enhanced by careful
pivoting can be effective in transforming the problem
to an unconstrained sparse-dense least squares problem that can be solved
with existing direct or iterative methods. We also present a number of solution variants
that employ an augmented system formulation, which can be attractive for solving
a sequence of related problems. Numerical experiments on problems coming
from practical applications are used throughout to demonstrate the
effectiveness of the different approaches.

\end{abstract}

\section{Introduction}
\label{intro}
Our interest lies in efficient and robust
methods for solving large-scale linear least squares  problems with linear equality constraints.
We assume that $A \in \Re^{m \times n}$ and  $C \in \Re^{p \times n}$, with
$m > n \gg p$. We further assume that
$A$ is large and sparse and $C$ represents a few, possibly dense, linear constraints.
Given $b \in \Re^{m}$ and $d \in \Re^{p}$, the least squares problem with equality
constraints (the LSE problem) is
\begin{equation}\label{eq:ls}
    \min_{x\in \Re^{n}} \left\| A x -b \right\|^2_2
\end{equation}
\begin{equation}\label{eq:constraints}
    \mbox{s.t.} \;\; C x = d.
\end{equation}
\noindent
A solution exists if and only if  (\ref{eq:constraints}) is consistent.
For simplicity, we assume that $C$ has full row rank (although the
proposed approaches can be made more general). In this case,
(\ref{eq:constraints}) is consistent for any $d$.
A solution to the LSE problem (\ref{eq:ls})--(\ref{eq:constraints}) is unique if and only if
${\cal N}(A) \cap {\cal N}(C) =\{0\}$, where
for any matrix $B$, ${\cal N}(B)$ denotes its null space. This is equivalent to the
extended matrix
\begin{equation}\label{eq:extended_matrix}
{\cal A} = \begin{pmatrix} A \\ C \end{pmatrix}
\end{equation}
having full column rank. In the case of non-uniqueness, there is a unique minimum-norm solution.

LSE problems arise in a variety of practical applications, including scattered data
approximation \cite{dast:13}, fitting curves to data \cite{fare:02}, surface fitting problems
\cite{pizi:07}, real-time signal processing, and control and communication leading to recursive problems \cite{zhli:07},
as well as when solving  nonlinear
least squares problems and least squares problems with inequality constraints.  For example,
in fitting curves to data, equality constraints may arise from the need to interpolate some
data or from a requirement for adjacent fitted curves to match with continuity of the curves
and possibly of some derivatives. Motivations for LSE problems together with solution strategies are summarized in
the research monographs \cite{bjor:96,bjor:15,laha:95}.

Classical approaches for solving LSE problems derive an equivalent unconstrained linear least squares (LS) problem
of lower dimension. There are two standard ways to perform this reduction: the null-space approach \cite{hala:69,laha:95} and the method of
direct elimination \cite{bjgo:67}, both of which, with suitable implementation,
offer good numerical stability. These methods, termed {\it constraint substitution} methods,
consider the constraints (\ref{eq:constraints})
as the primary data and substitute from them into the LS problem (\ref{eq:ls}).
The  former  performs a substitution
using a null-space basis of $C$ obtained from a QR factorization, while the latter
is based on substituting an expression for  selected solution components from the constraints into
(\ref{eq:ls}).  This can be done using either
a pivoted LU factorization \cite{bjgo:67} or a
QR factorization of $C$ \cite{bjor:96}.
Other  solution methods, which may be regarded as complementary
to the constraint substitution approaches, reverse the direction
of the substitution, substituting from the LS
problem into the constraints. This involves the use of an augmented
system and include  a Lagrange multiplier formulation \cite{heat:82},
updating procedures that force the constraints to be satisfied a posteriori \cite{bjor:84,bjor:96},
and a weighting approach \cite{baha:88,pore:69,vanl:85},

Solving large-scale LS problems is typically much harder than solving systems of linear
algebraic equations, in part because key issues such as
ill-conditioning or dense structures within an otherwise sparse problem can vary
significantly between different problem classes.
Consequently, we do not expect that there will be a single method that is optimal for all LSE problems, and having
a range of approaches available that target different problems is important.
Our main objective is to revisit classical solution strategies and to propose
new ideas and modifications that enable large-scale systems to be solved,
with an emphasis  first on the possibility that the constraints may be dense,
and second on requiring that the constraints be tightly satisfied.
In Sections~\ref{sec:null space} and \ref{sec:direct elimination},
we consider the null-space method and the
direct elimination approach, respectively. We review the methods
and show how they can be used for large-scale problems.
In Section~\ref{sec:aug system}, we  present complementary solution
approaches within an augmented system framework. This allows
us to treat the constraints and the least squares part of the problem using a single
extended system of equations or via a global updating scheme. Both direct and iterative
methods are discussed.

Much of the published literature related to LSE problems
lacks numerical results. For instance, Bj{\"o}rck \cite{bjor:84} remarks ``no attempt has yet been made to implement
the (general updating LSE) algorithm'', and as far as we are aware, attempts
remain absent. We assume this is because implementing the algorithms
is far from straightforward.
While it is not the intention here to offer a full general comparison of the different approaches,
throughout our study
we use numerical experiments on problems arising from real applications  to
highlight key features that may make a method attractive (or unsuitable) for particular problems
and to illustrate the effectiveness of the different approaches.
Our key findings and recommendations are summarized in Section~\ref{sec:conclusions}.

We end this introduction by describing our test environment.
The test matrices are taken from the SuiteSparse Matrix Collection \cite{dahu:2011}
and comprise a subset of those used by Gould and Scott in their study of numerical methods for
solving large-scale LS problems \cite{gosc:2017}.
If necessary, the matrix is transposed to give an overdetermined system.
Basic information on our test set is given in Table~\ref{T:test problems}.
\begin{table}[htbp]
\caption{Statistics for our test set.
$m$, $n$ and $nnz(\cal A)$ are, respectively, the row and column counts and the number of entries in the matrix ${\cal A}$
given by (\ref{eq:extended_matrix}). $dratio$ is
the ratio of the nonzero counts of the densest row to the sparsest row of ${\cal A}$. 
$^\dag$ indicates at least one column was removed to ensure there are no null columns in $A$.
}
\label{T:test problems}\vspace{3mm}
\begin{center}
\begin{tabular}{lrrrrrr} \hline
{Identifier} &
{$m$} &
{$n$} &
 {$nnz(\cal A)$} &
 {$dratio$} &
 {$\|x\|_2$~~~~} &
  {$\|r\|_2$~~~~} \\
\hline
lp\_fit2p &  13,525 & 3,000 & 50,284  & 3,000 &   1.689$\times 10^{1}$ &   1.105$\times 10^{2}$ \\
sc205-2r$^\dag$ & 62,423 & 35,212 & 123,237  & 1,602 &   8.758$\times 10^{1}$ &   2.039$\times 10^{2}$ \\
scagr7-2b$^\dag$ & 13,847 & 9,742 & 35,884 & 1,792 &   1.109$\times 10^{2}$ &   6.071$\times 10^{1}$ \\
scagr7-2r$^\dag$ & 46,679 & 32,846 & 120,140 & 6,048 &   1.821$\times 10^{2}$ &   1.133$\times 10^{2}$ \\
scrs8-2r$^\dag$ & 27,691 & 14,357 & 58,429 & 2,051 &   8.570$\times 10^{1}$ &   1.465$\times 10^{2}$ \\
sctap1-2b & 33,858 & 15,390 & 99,454  & 771 &   1.463$\times 10^{2}$ &   1.716$\times 10^{2}$ \\
sctap1-2r & 63,426 & 28,830 & 186,366  & 1,443 &   1.649$\times 10^{2}$ &   2.070$\times 10^{2}$ \\
south31 & 36,321 & 18,425 & 112,328 & 17,520 &   2.749$\times 10^{1}$ &   1.881$\times 10^{2}$ \\
testbig & 31,223 & 17,613 & 61,639 & 802 &   6.399$\times 10^{1}$ &   1.441$\times 10^{2}$ \\
\hline

deter3\_20   &  21,777 &   7,647 &  44,547 &     73 &   1.589$\times 10^{3}$ &   1.220$\times 10^{2}$ \\
deter3\_5    &  21,762 &   7,647 &  43,807 &     73 &   1.568$\times 10^{3}$ &   1.218$\times 10^{2}$ \\
fxm4\_6\_20   &  47,185 &  22,400 & 265,442 &     24 &   5.001$\times 10^{2}$ &   9.596$\times 10^{1}$ \\
fxm4\_6\_5    &  47,170 &  22,400 & 265,141 &     24 &   5.332$\times 10^{2}$ &   9.592$\times 10^{1}$ \\
gemat1\_20  & 10,595 & 4,929 & 47,369 &  22        & 3.170$\times 10^{4}$ &  8.595$\times 10^{1}$  \\
gemat1\_5  &  10,580 & 4,929 & 47,339 &  28        & 2.445$\times 10^{4}$ &  8.192$\times 10^{1}$  \\
stormg2-8\_20 &  11,322 &   4,393 &  28,553 &     21 &   2.829$\times 10^{1}$ &   7.970$\times 10^{1}$ \\
stormg2-8\_5 &  11,307 &   4,393 &  28,273 &     21 &   3.974$\times 10^{1}$ &   7.780$\times 10^{1}$ \\
\hline
\end{tabular}
\end{center}

\end{table}
The problems in the top half of the table contain rows that are identified as
dense by Algorithm~1 of \cite{sctu:2021a}  (with the density
parameter set to 0.05). These rows are taken to form
the constraint matrix $C$ and all other rows form $A$.
For the other problems, we form $A$ by removing the 20 densest rows of the
SuiteSparse matrix;
some or all of these rows are used to form
$C$ (and the rest are discarded).
Table~\ref{T:test problems} reports data for $p=5$ and 20 (denoted, for example, by deter\_5
and deter\_20, respectively).
Although the densest rows are not necessarily very dense, we make this choice because it corresponds
to the typical situation in which the constraints couple many of the solution components together.
For some of our test examples, splitting the supplied
matrix into a sparse part and a dense part results in the sparse part $A$  containing a small number of
null columns (at most 7 such columns for our test examples).
For the purpose of our experiments, we remove the corresponding columns from
the extended matrix (\ref{eq:extended_matrix})
(the data in Table~\ref{T:test problems} is for the modified problem).
In all our tests, we check that the norms of the computed
solution $x$ and least squares residual $r = b- Ax$ are consistent with the values
given in Table~\ref{T:test problems}.

In our experiments, we prescale the extended matrix $\mathcal{A}$
given by (\ref{eq:extended_matrix}) by normalizing each of
its columns. That is, we replace $\mathcal{A}$ by $\mathcal{AD}$, where $\mathcal{D}$
is the diagonal matrix with entries $\mathcal{D}_{ii}$
satisfying $\mathcal{D}_{ii} =1/ \|\mathcal{A}e_i\|_2 $ ($e_i$ denotes the $i$-th unit vector).
The entries of $\mathcal{AD}$ are at most one in absolute value.
The vectors $b$ and $d$ are set to be vectors of 1's (so that $\|b\|_2$ and
$\|d\|_2$ are $O(1)$).

For the substitution approaches described in Sections~\ref{sec:null space}
and \ref{sec:direct elimination}, we have developed prototype Fortran codes;
in Section~\ref{sec:aug system}, the augmented system methods
are implemented using the SuiteSparseQR package of
Davis~\cite{davi:11}
and Fortran software from the HSL mathematical software library
\cite{hsl:2018}.
The prototype codes are not optimised for efficiency and so computational times are not reported.
Developing library quality implementations is far from trivial and is outside the scope of the current study, which focuses
rather on determining which approaches are sufficiently promising for sophisticated
implementations to be considered in the future.

\medskip
\noindent
{\bf Notation}
All norms are 2-norms and in the rest of the paper,
to simply the notation, $\|.\|_2$ is denoted by $\|.\|$.
$I$ is used to denote the identity matrix of appropriate dimension.
The entries of any matrix $B$ are $(B)_{i,j}$
and its columns are denoted by $b_1,b_2, \ldots $.
The null space of $B$ is ${\cal N}(B)$ and $Z$ is used to denote a matrix
whose columns form a basis for the null space  (i.e.,
$Z$ satisfies $BZ = 0$).
Permutation matrices are denoted by $P$ (possibly with a subscript).
The  {\it normal
matrix} for (\ref{eq:ls}) is $H = A^TA$.

\section{The null-space approach}\label{sec:null space}
The null-space approach is a standard technique for solving least squares problems.
It is based on constructing a matrix
$Z \in  \Re^{n \times (n-p)}$ such that its columns form a basis for ${\cal N}(C)$. Any $x \in \Re^{n}$ satisfying the constraints can be written in the form
\begin{equation}\label{eq:particular}
    x=x_1 + Z x_2,
\end{equation}
where $x_1 \in \Re^{n}$ is a particular solution of the underdetermined system
$    C x_1 = d.$
The minimum norm solution can be obtained from the QR factorization of $C$,
that is, $CP_C = Q_C \begin{pmatrix}R_C & 0 \end{pmatrix}$, where 
the permutation $P_C \in \Re^{n \times n}$ represents
the pivoting, $R_C \in \Re^{p \times p}$ 
is an upper triangular matrix and $Q_C \in \Re^{p \times p}$
is an orthogonal matrix. $x_1$ is then given by
\[x_1 = P_C \begin{pmatrix} R_C^{-1}Q_C^Td \\ 0 \end{pmatrix}.\]
Substituting (\ref{eq:particular}) into (\ref{eq:ls}) gives the transformed LS problem
\begin{equation}\label{eq:transformed}
    \min_{x_2} \left\| A Z x_2 - (b -A x_1) \right\|^2.
\end{equation}
The method is summarized as Algorithm~\ref{alg:ns1}.
 \begin{algorithm}[htpb]
\caption{\label{alg:ns1}
Null-space method for solving the LSE problem (\ref{eq:ls})-(\ref{eq:constraints})
with $C$ of full row rank.}
\begin{algorithmic}[1]
\State Find $x_1 \in \Re^{n}$ such that $C x_1 = d$. 
\State Construct $Z \in \Re^{n \times (n-p)}$ of full column rank such that $CZ = 0$.
\State Solve the normal equations $Z^THZ x_2 = (AZ)^T(b -  A x_1)$ corresponding to (\ref{eq:transformed}) \Comment{Here $H = A^TA$}.
\State Set $x = x_1 + Zx_2$.
\end{algorithmic}
\end{algorithm}

In the 1970s, the null-space method was developed and discussed by a number of authors,
including in relation to quadratic programming
\cite{hala:69,laha:95,scst:79,stoe:71}.
These and subsequent contributions formulate the approach via the orthogonal null-space basis
obtained, for example, from the QR factorization of $C^T$ given by
\[C^T = PQ \begin{pmatrix} R \\ 0 \end{pmatrix}, \]
where $P \in \Re^{n \times n}$ is a column permutation of $C$ and $Q \in \Re^{n \times n}$
is an orthogonal matrix. $Z$ is equal to the last $(n-p)$ columns of $PQ$
and consequently is dense. Note that although it is possible to store $Q$ implicitly 
using, for example, Householder transformations, the memory
demands and implied operation counts are generally too high.
Our interest is in large LS problems and therefore
it  may not be practical to solve
the $(n-p) \times (n-p)$ system in Step 3 if $Z$ is dense.
To make the approach feasible for large problems
we can exploit our recent work \cite{sctu:2021b}
on constructing sparse null-space bases of ``wide'' matrices
such as $C$ that have many more columns than rows
and may include some dense rows.

Scott and T\accent23uma~\cite{sctu:2021b} propose a number of ways to construct sparse $Z$.
In our experiments, we employ Algorithm~3 from Section 3 of \cite{sctu:2021b}.
This computes each column of $Z$ independently using its own QR factorization
that incorporates a threshold pivoting strategy. 
When selecting candidate pivots, the aim is to 
combine using a threshold parameter $\theta \in [0,1]$ 
  to ensure stability
with locality (that is, column interchanges are between 
columns that are as close as possible) to maintain sparsity in
$Z$. Small values of $\theta$ lead to $Z$ having a narrow bandwidth and (assuming
the normal matrix $H$ is sparse)
to $Z^THZ$ being sparse but this is potentially at the  expense of numerical stability.

\begin{figure}[htbp]
\begin{center}
\includegraphics[height=4.5cm]{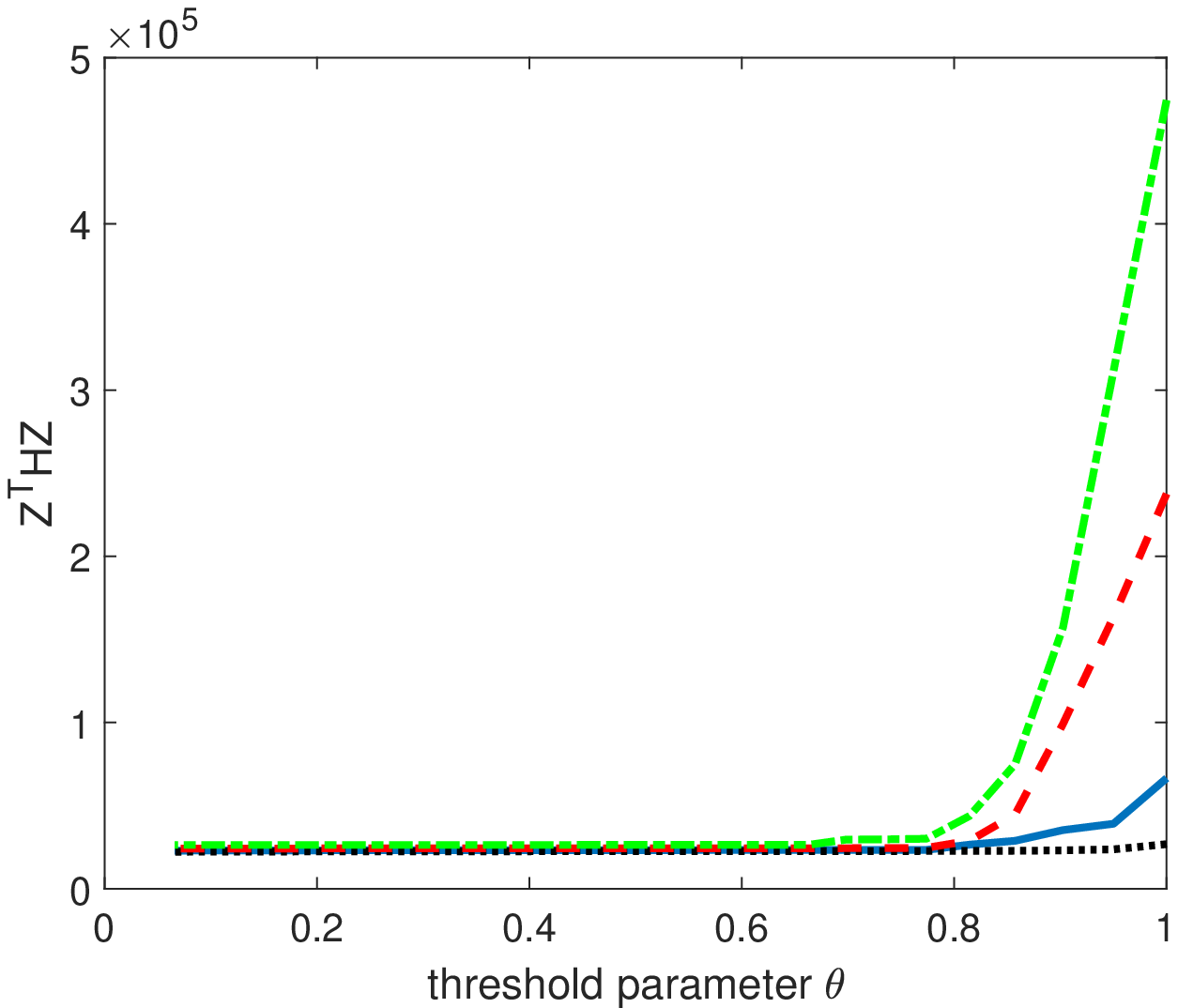} \hskip1.cm
\includegraphics[height=4.5cm]{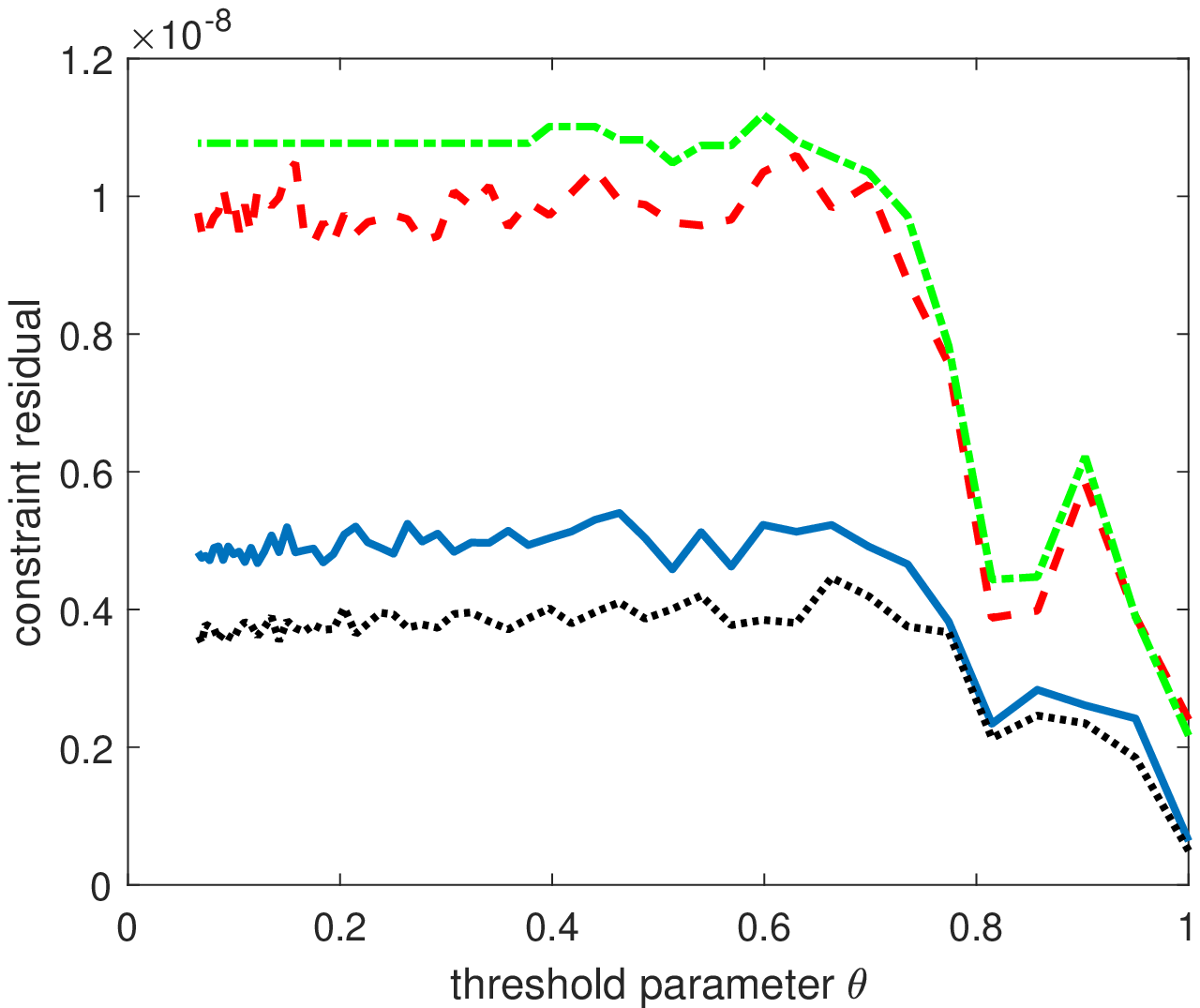} \\
\vspace{0.65cm}
\includegraphics[height=4.35cm]{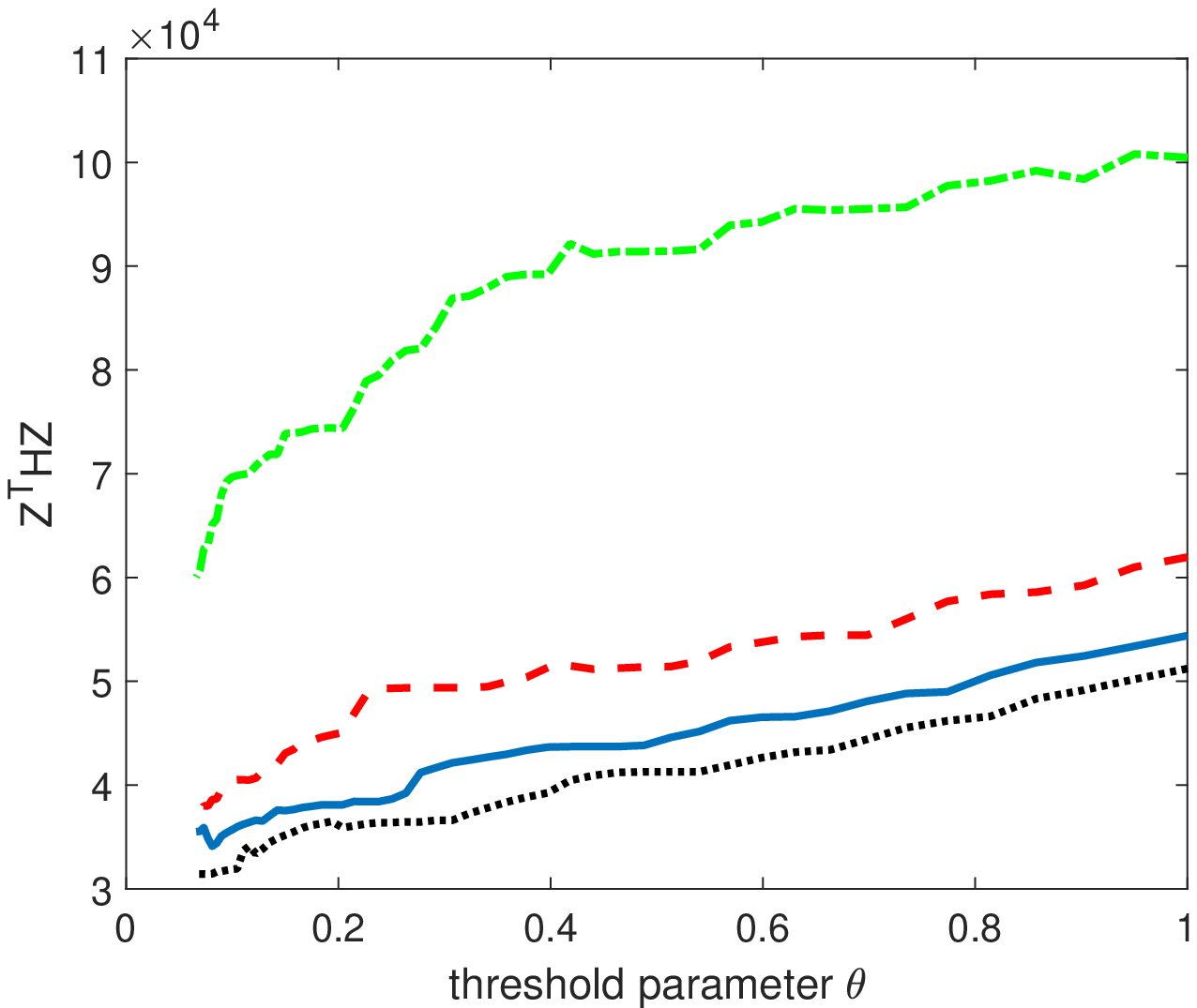} \hskip1.cm
\includegraphics[height=4.25cm]{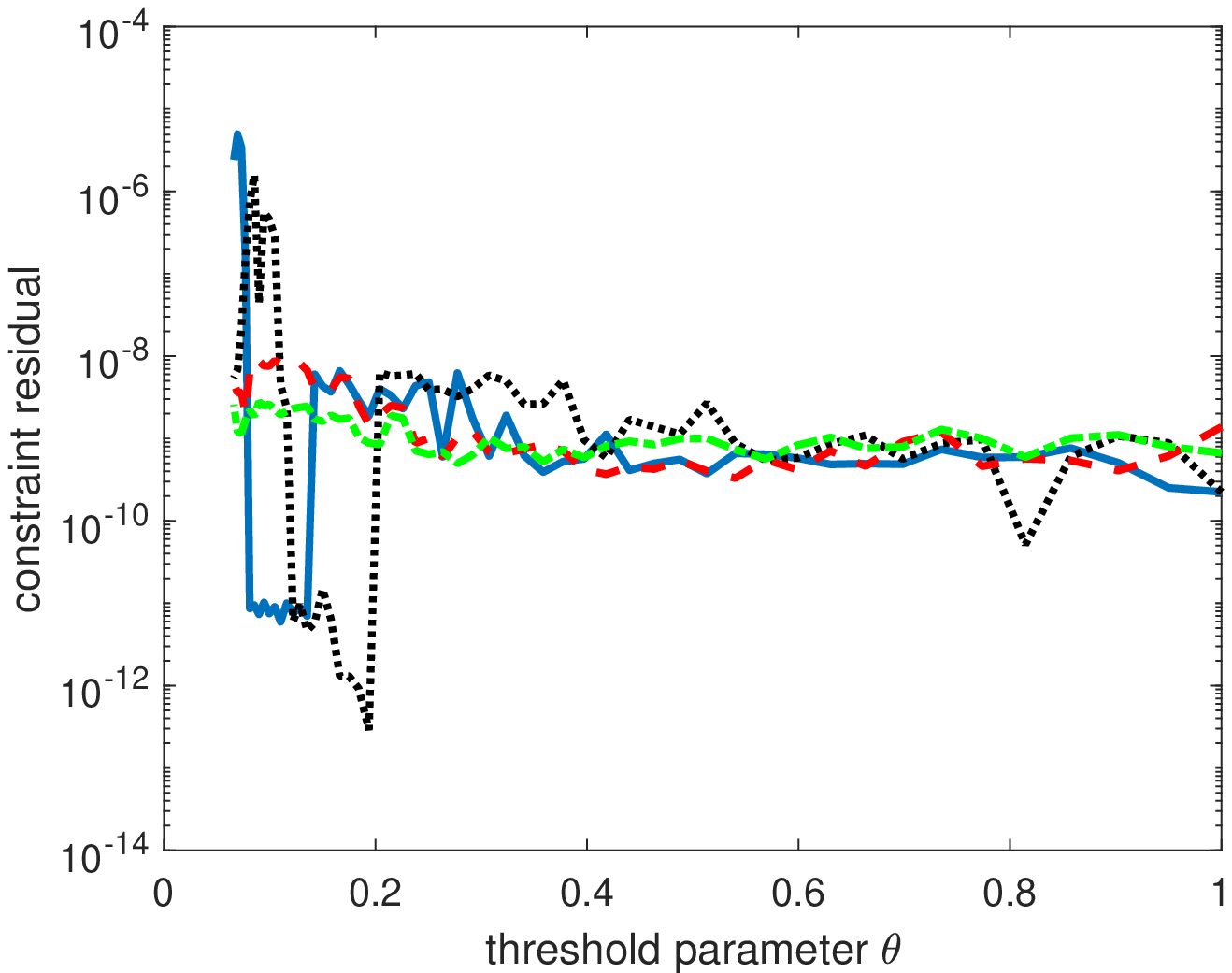}
\end{center}
\caption{The number of entries in $Z^THZ$ (left) and the constraints residual
$\left\| r_c\right\|$ (right) for problem deter3 (top) and gemat1 (bottom)
as the threshold pivoting parameter $\theta$ used in the computation of the null-space basis increases from 0.3 to 1.
The four curves correspond to $p=2$ (black dotted line), 5 (blue full line),
10 (red dashed line) and 20 (green dash-dotted line).
\label{fig:deter3 and gemat1}
}
\end{figure}

\begin{table}[htbp]
\caption{The  density of $Z^THZ$ (that is, $nnz(Z^THZ)/(n-p)^2$)  and
constraint residual $\left\| r_c \right\|$
for two values of the threshold pivoting parameter $\theta$  used
in the computation of the null-space basis.
$\ddag$ indicates insufficient memory for {\tt HSL\_MA87}.
\label{fig:NS_table}
}
\begin{center}
\begin{tabular}{lr cccc}   \hline
\multicolumn{1}{c}{}  &
\multicolumn{1}{c}{}  &
    \multicolumn{2}{c}{$\theta=1$} &

    \multicolumn{2}{c}{$\theta=0.1$}  \\
  \hline
Identifier & $p$  &  Density   &   $\left\| r_c \right\|$ & Density  &   $\left\| r_c \right\|$  \\
\hline
lp\_fit2p   &  25 & 0.47 & 4.14$\times 10^{-5}$ & 0.11 & 3.07$\times 10^{-5}$ \\
sc205-2r & 8 & 0.03 & 5.25$\times 10^{-8}$ & ~~~0.0002 & ~7.58$\times 10^{-11}$ \\
scagr7-2b & 7 & 0.03 & 1.27$\times 10^{-8}$ & ~~~0.0007 & 2.79$\times 10^{-8}$ \\
scrs8-2c & 22 & 0.31 & ~3.60$\times 10^{-11}$& 0.23 & ~2.02$\times 10^{-11}$ \\
sctap1-2b  &  34 & 0.05 & 3.67$\times 10^{-6}$ & ~0.002 & 4.37$\times 10^{-7}$\\
sctap1-2r  &  34 & 0.05 & 2.76$\times 10^{-3}$ & 0.02 & 1.83$\times 10^{-4}$ \\
south31  & 5 & 0.20 & $\ddag$ & 0.02 & 3.26$\times 10^{-7}$ \\
testbig & 8 & 0.03 & ~2.53$\times 10^{-11}$ & ~~~0.0002 & ~2.90$\times 10^{-11}$ \\
\hline
deter3\_20 &  20 &  0.008 & 2.58$\times 10^{-9}$ &  ~0.0004 &  1.09$\times 10^{-8}$ \\
deter3\_5 & 5 & 0.001 & 6.39$\times 10^{-10}$ & ~0.0004 & 4.89$\times 10^{-9}$ \\



fxm4\_6\_20 & 20 &  ~0.0006 &  5.43$\times 10^{-6}$ &  ~0.0006 & 6.67$\times 10^{-7}$ \\
fxm4\_6\_5 & 5 & ~0.0005 & ~7.80$\times 10^{-11}$ & ~0.0005 & ~1.10$\times 10^{-11}$ \\

gemat1\_20 & 20 & 0.004  &  1.10$\times 10^{-9}$&  0.003 &  2.40$\times 10^{-9}$ \\
gemat1\_5 & 5 & 0.002  & ~2.24$\times 10^{-10}$& 0.001 & ~1.00$\times 10^{-11}$ \\




stormg2-8\_20 & 20 & 0.003 &  7.44$\times 10^{-9}$ &  0.002 &  7.23$\times 10^{-9}$\\
stormg2-8\_5 & 5 & 0.002 & ~1.13$\times 10^{-10}$ &  0.002 & ~8.16$\times 10^{-11}$\\


 \hline
\end{tabular}
\end{center}
\end{table}
Our first results are for problems deter3 and gemat1.
As discussed in the Introduction,
we form the constraint matrix $C$ by taking the $p = 2$, 5, 10, 20  densest rows
of $\mathcal A$. The sparse block
$A$ is the same for each case.
In Figure~\ref{fig:deter3 and gemat1}, we plot the number of entries $nnz(Z^THZ)$ in $Z^THZ$ and the
norm of the constraints residual $\left\| r_c\right\| =  \left\| d - C x \right\|$.
As expected, $nnz(Z^THZ)$ increases with $\theta$, and this increase grows with $p$.
This is illustrated further by the results in  Table~\ref{fig:NS_table}.
We can see that, independently of the choice of $\theta$, for some problems
(including lp\_fit2p and sctap1-2r) the constraints are not
tightly satisfied.  This demonstrates an inherent limitation
of the null-space approach of \cite{sctu:2021b} 
that focuses on constructing the columns of $Z$ so as to keep
$Z^THZ$ sparse but does not result in $Z$ having orthogonal columns. 

The linear system in Step 3 of Algorithm~\ref{alg:ns1}
is symmetric positive definite. In the above experiments, we employ
the sparse direct solver {\tt HSL\_MA87} \cite{hors:2010} (combined with an approximate minimum
degree ordering).
However, for large problems, the memory demands mean it may not be possible to use a direct method;
this is illustrated by problem south31 with $\theta = 1$.
If a preconditioned iterative solver is used instead, not only are the solver memory requirements
much less but explicitly forming the potentially ill-conditioned
normal matrix $H$ can be avoided and because $Z$ only needs to be applied implicitly,
the need for sparsity can potentially be relaxed.
Currently, finding a good preconditioner for use in this case
remains an open problem \cite{naso:96}.

If a sequence of LSE problems is to be solved with the
same set of constraints but different $A$, the null-space basis can be reused,
substantially reducing the work required. But if the constraints
are changed, then  $Z$ will also change. In \cite{sctu:2021b}, we present a strategy that allows
$Z$ to be updated when a row (or block  of rows) is added to $C$.

\section{The method of direct elimination}\label{sec:direct elimination}

The second method we look at is direct elimination \cite{laha:95}.
The basic idea is to express the dependency of $p$ selected  components of the vector $x$
on the remaining $n-p$ components and to substitute this  into the
 LS problem (\ref{eq:ls}).  Here we propose how to choose the $p$ components so as to retain
sparsity in the transformed problem.

Consider the constraints (\ref{eq:constraints}). The method starts by permuting and splitting
the solution components  as follows:
\begin{equation*}
Cx = CP_c y =
    \begin{pmatrix}
C_1 & C_2
    \end{pmatrix}
    \begin{pmatrix}
y_{1} \\ y_{2}
    \end{pmatrix} = d,
\end{equation*}
where $P_c\in \Re^{n \times n} $ is a permutation matrix
chosen so that $C_1 \in \Re^{p \times p}$ is nonsingular. Let $AP_c = \begin{pmatrix}
A_1 & A_2     \end{pmatrix}$
be a conformal partitioning of $AP_c$. Substituting
the  expression
\begin{equation}\label{eq:y1}
y_1 = C_1^{-1}(d - C_2y_2)
\end{equation}
into (\ref{eq:ls}) gives the transformed LS problem
\begin{equation}\label{eq:ls2}
    \min_{y_2} \left\| A_T y_2 -(b - A_1C_1^{-1}d)\right\|^2_2,
\end{equation}
with the transformed matrix
\begin{equation}\label{eq:tranformed mx}
A_T = A_2 - A_1C_1^{-1}C_2 \in \Re^{m \times (n-p)}.
\end{equation}
Note that if $C_1$ is irreducible, the transformation combines all the rows of $C_2$.
If $C$ is composed of dense rows then $A_T$ has  more dense rows than $A$.  We thus
seek to add as few row patterns as possible  replicating the (possibly) dense pattern of $C$ within $A_T$.
If both $A$ and $C$ are sparse, the substitution leads to a sparse LS problem.
We have the following straightforward result.
\medskip
\begin{lemma} \label{lm:de}
Let $A \in \Re^{m \times n}$ be sparse.
Let $m  >n > p$ and assume a conformal column splitting
induced by the
permutation $P_c$ is such that
$CP_c = \begin{pmatrix}C_1 & C_2\end{pmatrix}$  and $AP_c = \begin{pmatrix}A_1 & A_2\end{pmatrix}$
with $C_1 \in \Re^{p \times p}$ nonsingular and $A_1 \in \Re^{m \times p}$.
Define the index set
\[
 Occupied =\{i \, | \, (A_1)_{i,k} \neq 0 \mbox{  for some $k$, } 1 \le k \le p \}.
\]
Then the number of dense rows in $A_T$ given by (\ref{eq:tranformed mx}) is at most
the number of entries in $Occupied$.
\end{lemma}

\noindent
\begin{proof}
The result follows directly from the transformation. Assuming the rows of
$C_1^{-1}C_2$  are dense, the substitution step (\ref{eq:y1}) of the direct elimination
implies a dense row $k$ in
$A_T$ only if there is a nonzero in the $k$-th row of $A_1$.
\end{proof}

\medskip
A simple  example is given in Figure~\ref{fig:de_example}. Here  we
ignore cancellation of nonzeros during arithmetic operations.
We see that the pattern of $A_T$ satisfies Lemma~\ref{lm:de}.
Note that, although in this example $C_1^{-1}C_2$ is shown
as dense, it need not be fully dense and
the number of entries in $Occupied$ represents an upper bound on the number of dense rows in  $A_T$.
\begin{figure}[htbp]
\begin{center}
$$
     \begin{pmatrix} * \\
&& * &&& * \\
 & * &&&&&\\
&* \\
&&* & *\\
&* \\
&&&* &  & * \\
&&* & * &&* & *\\
*& &&& *
 \end{pmatrix}-
    \begin{pmatrix} * \\
&& * \\
\\
 \\
\\
&* \\
\\
\\
*
 \end{pmatrix}
    \begin{pmatrix}
 * & * & * & * & * & * & *  \\
 * & * & * & * & * & * & *  \\
 * & * & * & * & * & * & *  \\
 \end{pmatrix}
 \longrightarrow
\begin{pmatrix} * &* &* &*  &* &* &* \\
 * &* &* &* &* &* &*   \\
&*\\
 * &* &* &* &* &* &*  \\
&&* & *\\
&* \\
&&&* &  & * \\
&&* & * &&* & *\\
 * &* &* &* &* &* &*
 \end{pmatrix}
$$
 \caption{\label{fig:de_example} Example of the transformation in the direct elimination approach.
Here $m = 9$, $p=3$, $n=7$. The depicted matrices (from the left)
represent the transformation  $A_T = A_2 - A_1(C_1^{-1}C_2)$. The matrix
$C_1^{-1}C_2 \in \Re^{p \times n} $ is given as fully dense.}
\end{center}
\end{figure}

Lemma~\ref{lm:de} implies that the LSE problem is transformed to a LS problem (\ref{eq:ls2})
that has some dense rows, which we refer to as a {\em sparse-dense} LS problem. Consequently,
existing methods for  sparse-dense LS problems can be used, including
those recently proposed in \cite{sctu:2017b,sctu:2018a,sctu:2020a}
(see also the recent direct LS solver {\tt HSL\_MA85} from the HSL
library).
A straightforward algorithmic implication of the lemma is that the permuting and
splitting of $C$ cannot be separated from considering
the sparsity pattern of $A$ because the splitting also
determines $A_1$ and $A_2$. Thus we want to
permute the columns of $C$ to allow a sufficiently
well-conditioned factorization of $C_1$ while limiting the number of entries in $Occupied$
and hence the number of dense rows in  $A_T$.
The approach outlined in Algorithm~\ref{alg:de1} is one way to achieve this.
There is an important difference between the
pivoting used in Algorithm~3 of \cite{sctu:2021b} (which we used in the previous section) and 
that of Algorithm~\ref{alg:de1} below.
The former  incorporates standard threshold pivoting into 
a QR factorization algorithm that aims to compute the 
columns of $Z$ with very few nonzero entries so that the transformed normal
matrix $Z^THZ$ is as sparse as possible. Its  threshold 
parameter $\theta$ controls the 
size of the pivots and thus the stability of the factorization. 
The threshold parameter $\tau \in (0,1]$ 
used in Algorithm~\ref{alg:de1} also guarantees the pivots
in the QR factorization of  $C$ are not small but 
the selection of the candidate pivots is very different.
In this case, the choice
is based on limiting the fill-in in the transformed matrix $A_T$.
Here it is the set of rows held in $Occupied$ that
potentially fill-in in $A_T$ that plays a crucial role.
The use of different notation for the threshold parameters
emphasises the difference between the two QR-based approaches. 


Observe that the pivoting strategy
in Algorithm~\ref{alg:de1} considers $C$ and $A$ simultaneously and
will not select a column as the pivot column if this column in $A$ is dense (as
it would lead to $A_T$ being dense).
While we do not discuss the implementation details, we remark that
care is needed to ensure efficiency. For example,
QR with pivoting for a wide matrix is relatively cheap but it may be necessary
to store the squares of the column norms using a heap, which is why we emphasize their role in the algorithm
by using the explicit notation $w_i$ for these norms.
\begin{algorithm}[htpb]
\caption{\label{alg:de1}
Assume $C = ( c_1, \ldots, c_n )\in \Re^{p \times n}$ ($p < n$) has full row rank.
Compute $C_1 \in \Re^{p \times p}$  and the
column permutation $P_c\in \Re^{n \times n}$ for the direct elimination method for solving the LSE problem (\ref{eq:ls})--(\ref{eq:constraints}).
$P_c$ is determined by a QR factorization with threshold pivoting; 
$\tau \in (0,1]$ is the threshold pivoting parameter.}
\begin{algorithmic}[1]
\State Initialise: $Occupied = \emptyset$, $S = \emptyset$, and $w_i = \|c_i\|^2$, $i=1, \ldots, n$. Define $E_n = \{1,2, \ldots , n\}$.
\State {\bf for} $l = 1,...,p$
\State \indent Find $i_{max} = \argmax_{i \in E_n} \{w_i \ | \ i\in E_n \setminus S\}$.
\State \indent Define $E_{\tau} = \{i \in E_n \setminus S \ | \ w_i \ge \tau w_{imax}\}$.
\State \indent Find $i_l  \in E_n \setminus S$ such that
$i_l=\argmin_{k \in E_{\tau} }|\{j \, | \, (A)_{j,k} \neq 0 \} \setminus Occupied|$.
\State \indent For $j \in E_n \setminus S$ set $c_j \leftarrow c_j - q_{i_l}^Tc_j \, q_{i_l}$, where $q_{i_l} = c_{i_l}/\|c_{i_l}\|$.
\State \indent For $j \in E_n \setminus S$ set $w_j \leftarrow w_j - (q_{i_l}^Tc_j)^2$.
\State \indent Update $S \leftarrow S \cup \{i_l\}$.
\State \indent Update $Occupied \leftarrow Occupied \cup \{j \, | \, (A)_{j,i_l} \neq 0 \}$.
\State {\bf end for}
\State Set $P_c$ to permute the columns of $C$ with indices in $S$ to obtain $C_1$.
\end{algorithmic}
\end{algorithm}

\begin{figure}
\begin{center}
\includegraphics[height=4.3cm]{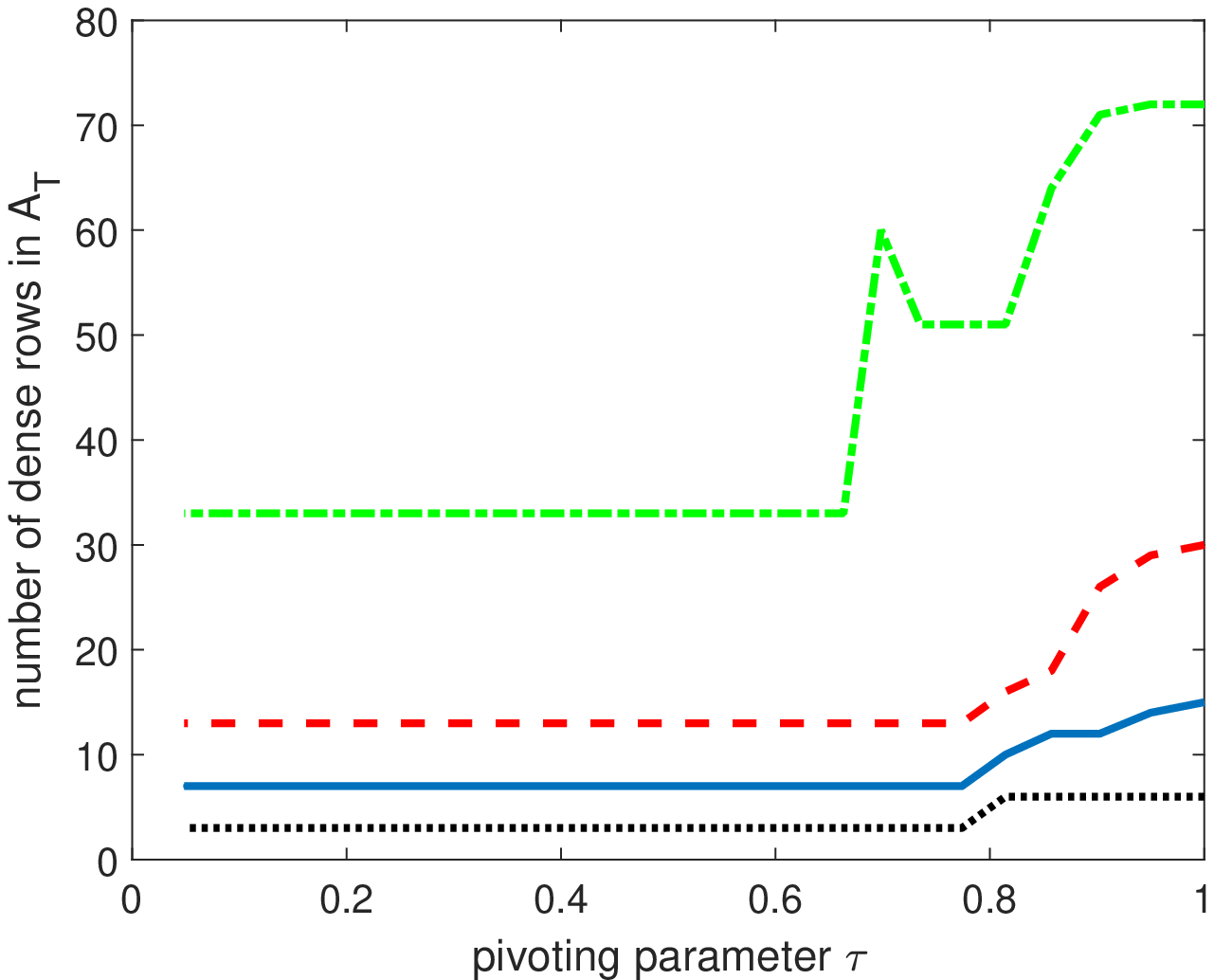} \hskip1.cm
\includegraphics[height=4.3cm]{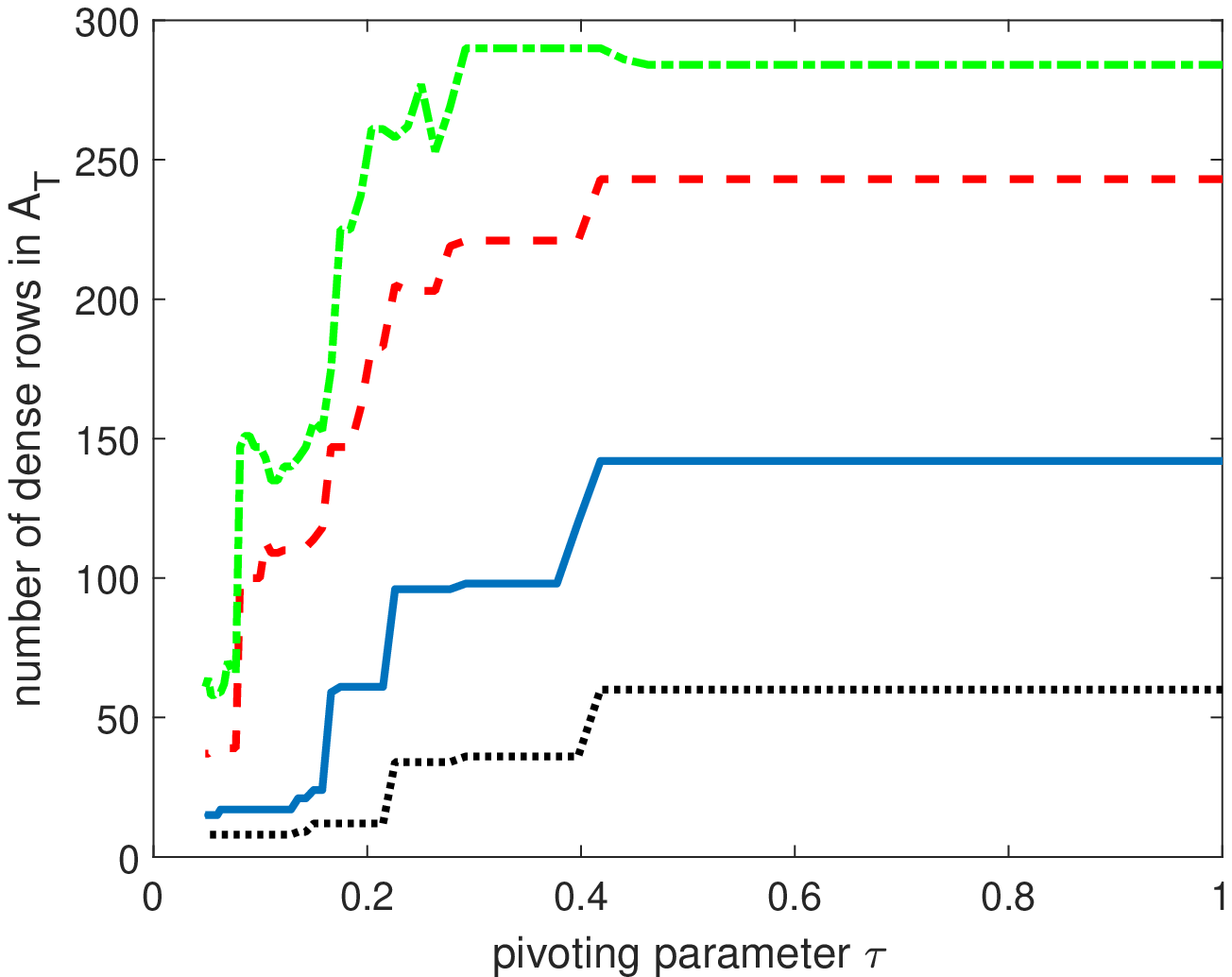}
\end{center}
\caption{The number of dense rows in the transformed matrix $A_T$
as the parameter $\tau$ increases from 0.05 to 1 for problems deter3 (left) and gemat1 (right).
The four curves correspond to
$p=2$ (black dotted line), 5 (blue full line), 10 (red dashed line) and 20 (green dash-dotted line).
}
\label{fig:deter3_DE}
\end{figure}
\begin{table}[htbp]
\caption{The number $ndense$ of dense rows in $A_T$ and norm of the constraints
residual $\left\| r_c \right\|$ for
two values of the pivoting parameter $\tau$.
\label{fig:DE_table}
}
\begin{center}
\begin{tabular}{ll rcrc}   \hline
\multicolumn{1}{c}{}  &
\multicolumn{1}{c}{}  &
    \multicolumn{2}{c}{$\tau=1$} &
    \multicolumn{2}{c}{$\tau=0.1$}  \\
  \hline
Identifier & $p$  &  $ndense$   &   $\left\| r_c \right\|$ & $ndense$   &   $\left\| r_c \right\|$  \\
\hline
lp\_fit2p  &  25 & 115 & 8.12$\times 10^{-12}$ & 100 &  6.77$\times 10^{-11}$ \\
sc205-2r & 8 & 9 & 3.32$\times 10^{-14}$ & 7 & 6.13$\times 10^{-14}$ \\
scagr7-2b & 7 & 14 & 4.03$\times 10^{-13}$ & 6 & 1.48$\times 10^{-13}$ \\
scrs8-2c & 22 & 16 & 6.80$\times 10^{-14}$  & 16 & 1.25$\times 10^{-13}$ \\
sctap1-2b  &  34 & 72 & 7.60$\times 10^{-13}$ & 63 & 1.31$\times 10^{-12}$\\
sctap1-2r  &  34 & 66 & 1.90$\times 10^{-13}$ & 57 & 1.40$\times 10^{-13}$\\
south31  & 5 & 20 & 2.57$\times 10^{-15}$ & 16 & 3.45$\times 10^{-14}$ \\
testbig & 8 & 9 &4.09$\times 10^{-15}$& 8 & 1.06$\times 10^{-14}$ \\
\hline
deter3\_20 & 20 & 72 &  8.02$\times 10^{-14}$ &  33 &  5.75$\times 10^{-14}$ \\
deter3\_5  & 5 & 15 &  8.02$\times 10^{-14}$ & 7 & 1.53$\times 10^{-13}$ \\



fxm4\_6\_20 & 20 & 113  &  2.07$\times 10^{-14}$ & 80 &  5.70$\times 10^{-14}$ \\
fxm4\_6\_5 & 5 & 50  & 7.81$\times 10^{-16}$ & 14 & 9.13$\times 10^{-16}$ \\

gemat1\_20 & 20 &  284 &  8.95$\times 10^{-15}$ & 147 &  1.64$\times 10^{-14}$  \\
gemat1\_5 & 5 & 142 & 9.77$\times 10^{-14}$ & 17 & 4.61$\times 10^{-14}$  \\




stormg2-8\_20 &  20 &  136 & 5.30$\times 10^{-14}$&  94 &  2.29$\times 10^{-15}$\\
stormg2-8\_5 & 5   & 61 & 3.28$\times 10^{-15}$& 35 & 1.50$\times 10^{-14}$\\


 \hline
\end{tabular}
\end{center}
\end{table}
The effects of increasing the  pivoting parameter $\tau$ on the number of dense rows in $A_T$
are illustrated in Figure~\ref{fig:deter3_DE} for problems deter3 and gemat1;
results for the full test set are given in Table~\ref{fig:DE_table}.
The dense rows of the transformed matrix $A_T$ are determined using Algorithm~1 of \cite{sctu:2021a}
and to solve the transformed LS problem (\ref{eq:ls2}) we use the sparse-dense preconditioned iterative approach of \cite{sctu:2017b}.
This computes
a Cholesky factorization of the normal matrix corresponding to the
sparse part of $A_T$ and uses it as a preconditioner within a conjugate gradient (CG) method; the CG convergence tolerance
that measures relative decrease of the transformed residual 
$||A_T^Tr||_2/||r||_2$ is set to $10^{-11}$.
For the problems in the top half of the table for which the rows of
$C$ are much denser than those of $A$ (recall Table~\ref{T:test problems}), reducing
$\tau$ leads to only a small reduction in the number $ndense$ of
dense rows in $A_T$. However, when the constraints are not dense (the problems in the lower
half of the table), $ndense$ can be significantly decreased by choosing $\tau < 1$,
although if $\tau$ is too small, 
the matrix $C_1$ computed by Algorithm~\ref{alg:de1} can
become highly ill-conditioned and $A_T$ close to being singular.
In our experiments we occasionally observed this for $\tau < 10^{-5}$.

By comparing the pairs of problems in the lower half of the table (such
as deter3\_5 and deter3\_20) and
considering the plots in Figure~\ref{fig:deter3_DE},
we see that increasing the number $p$ of constraints can lead to a sharp increase in
$ndense$ (even if these constraints are relatively sparse), which can result in the transformed problem being hard to solve.
The constraints are
very well satisfied in all the successful test cases, making this an attractive
approach if a good sparse-dense LS solver is available and
the number of dense rows in the transformed problem is not too large.
Furthermore, it can be used, without modification, if the matrix $A$
contains a (small) number of dense rows.
However, for a sequence of problems, if $A$ and/or $C$ changes
then, because direct elimination couples the two
matrices, the computation must be completely restarted.

\section{Approaches described via augmented systems}
\label{sec:aug system}
We now focus on complementary approaches that are based on substitution from the unconstrained least squares
problem into the constraints.
A useful way to describe this is via the augmented (or saddle-point) system
\begin{equation}
  \label{eq:saddle}
  \begin{pmatrix}
    H & C^T \\ C & 0
 \end{pmatrix}
  \begin{pmatrix}
    x \\ \lambda
  \end{pmatrix}
  =
  \begin{pmatrix}
    A^T b \\ d
  \end{pmatrix},  \qquad H = A^TA.
\end{equation}
Here $\lambda \in \Re^{p}$ is a vector
of additional variables that are often called {\em Lagrange multipliers} \cite{govl:13,heat:82}.
The solution $x$ of (\ref{eq:saddle}) solves the LSE problem.
Using (\ref{eq:saddle}) can be particularly useful if the constraint matrix $C$ is dense and $p$ is small. As we shall see in the following
discussions, this is because
the work involved in the proposed algorithms that depends upon $p$
is effectively independent of the density of $C$.
Observe that because (\ref{eq:saddle}) has a zero $(2,2)$ block,
the augmented system can be also used to give an alternative
derivation of the null-space approach of Algorithm~\ref{alg:ns1}.
For if $Z$ is such that $CZ = 0$ and $x_1$ is a particular solution of
the second equation of (\ref{eq:saddle}) so that $Cx_1 = d$
(steps 1 and 2 of Algorithm~\ref{alg:ns1}),
then if $x = x_1 + \hat x$, (\ref{eq:saddle}) becomes
\begin{equation*}
  \begin{pmatrix}
    H & C^T \\ C & 0
 \end{pmatrix}
  \begin{pmatrix}
    \hat x \\ \lambda
  \end{pmatrix}
  =
  \begin{pmatrix}
    A^T (b- Ax_1) \\ 0
  \end{pmatrix}.
\end{equation*}
The second equation in this system is equivalent to finding $x_2$ such that $\hat x = Zx_2$. Substituting
this into the first equation we have
\[ HZx_2 + C^T \lambda = A^T (b- Ax_1).\]
Hence
\[Z^THZx_2 = (AZ)^T (b- Ax_1) \]
(see steps 3 and 4 of Algorithm~\ref{alg:ns1}).

\subsection{Direct use of Lagrange multipliers}
\label{sec:Direct use of Lagrange multipliers}
Algorithm~\ref{alg:lagrange} presents a straightforward updating scheme for solving the LSE problem using
Lagrange multipliers and (\ref{eq:saddle}).
Any appropriate direct or iterative method can be used for Step 1, which is usually the most expensive part of the computation. There is
no dependence on $C$ so the solution $y$ does not need to be recomputed when $C$ changes.
The method used to solve the system with a block of $p$ right-hand sides in Step 2 can be chosen
to exploit Step 1. For example, a sparse Cholesky factorization of $H$ may be computed in Step 1 and
the factors reused in Step 2.
Using existing sparse LS solvers (and a dense linear solver
for the $p \times p$ at Step 5), Algorithm~\ref{alg:lagrange} is
straightforward to implement and, from Step 6, the
solution $y$ of the unconstrained LS problem can be compared with that of the LSE.
\begin{algorithm}[htpb]
\caption{\label{alg:lagrange}
Straightforward updating approach based on Lagrange multipliers for solving the LSE problem (\ref{eq:ls})-(\ref{eq:constraints})
with $C$ having full row rank}
\medskip
\begin{algorithmic}[1]
\State  Solve  the sparse unconstrained LS problem $   \min_y \left\| A y -b \right\|^2$
\State Solve
 $H J = -C^T$ for $J \in \Re^{n \times p}$.
\State Set  $Y=CJ$.
\State Solve $Y\lambda = d -Cy$ for $\lambda$. \Comment{Note that $Y \in \Re^{p \times p}$ is symmetric negative definite.}
 \State Set $ x = y + J\lambda$.
\end{algorithmic}
\end{algorithm}

As discussed by Golub~\cite{golu:65} and Heath~\cite{heat:82}, a
numerically superior direct method that avoids both forming the potentially ill-conditioned normal matrix
$H$ and computing the multipliers $\lambda$ can be derived using the QR factorization of $A$.
Following \cite{sctu:2020a}, we obtain Algorithm~\ref{alg:updating_LSE}. Here $P$ is a permutation matrix
chosen to ensure sparsity of the R factor. Note that, unless $b$ (and hence $f$)
changes, the Q factor need not be retained and the R factor can be reused
if the constraints change but $A$ is fixed. 
\begin{algorithm}[htpb]
\caption{\label{alg:updating_LSE}
QR algorithm with updating for solving the LSE problem (\ref{eq:ls})--(\ref{eq:constraints})
with $C$ having full row rank}
\medskip
\begin{algorithmic}[1]
\State Compute the QR factorization
$\begin{pmatrix} AP & b  \end{pmatrix} = Q\begin{pmatrix} R & f \\ 0 & g \end{pmatrix}$ using a sparse QR solver.
\State Solve $R P^T y = f$ for $y$.
\State Solve $P R^TK^T = C^T$ for $K^T \in \Re^{n \times p}$.
\State Compute the minimum-norm solution of $Ku = d - Cy$.
\State Solve $R P^T z = u$ for $z$.
\State Set $x = y + z$.
\end{algorithmic}
\end{algorithm}
\begin{table}[htbp]
\caption{Norm of the constraint residuals $\|r_c\|$ for QR with updating (Algorithm~\ref{alg:updating_LSE}).
\label{T:res updating}
}
\begin{center}
\begin{tabular}{l r | l  r| l  r}   \hline
\multicolumn{1}{c} {Identifier} &
\multicolumn{1}{c|}{$\|r_c\|$} &
\multicolumn{1}{c} {Identifier} &
\multicolumn{1}{c|}{$\|r_c\|$} &
\multicolumn{1}{c} {Identifier} &
\multicolumn{1}{c}{$\|r_c\|$} \\
  \hline
lp\_fit2p  &      4.485$\times 10^{-11}$  & sctap1-2b  &      4.422$\times 10^{-11}$  & deter3\_20   & 1.264$\times 10^{-12}$\\
sc205-2r   &      4.299$\times 10^{-10}$  & sctap1-2r  &      7.624$\times 10^{-11}$  & fxm4\_6\_20  & 8.493$\times 10^{-14}$ \\
scagr7-2b  &      1.364$\times 10^{-11}$  & south31    &      4.502$\times 10^{-13}$  & gemat1\_20 & 2.943$\times 10^{-12}$ \\
scagr7-2r  &      2.177$\times 10^{-11}$  & testbig    &      8.427$\times 10^{-11}$  & stormg2-8\_20 & 2.437$\times 10^{-12}$ \\
 scrs8-2r  &      8.634$\times 10^{-11}$  &            &                              &  & \\
 \hline
\end{tabular}
\end{center}
\end{table}
Results for Algorithm~\ref{alg:updating_LSE}  presented in Table~\ref{T:res updating}
confirm that the computed solution is such that the norm of
the constraints residual $\|r_c\|$ is small.
We omit results for problems such as deter\_5 that have $p=5$
constraints because they are similar (with $\|r_c\|$  typically smaller than for the corresponding problems with $p=20$).

\subsection{An extended augmented system approach}

An equivalent formulation of (\ref{eq:saddle}) is given by the 3-block saddle-point system (the first order
optimality conditions)
\[
{\mathcal A}_{aug} y= b_{aug},
\]
where
\begin{equation} \label{eq:LSE aug}
{\mathcal A}_{aug} =
\begin{pmatrix} I & 0 & A \\ 0 & 0 & C \\ A^T & C^T & 0 \end{pmatrix}, \qquad
y = \begin{pmatrix} r_s \\ \lambda \\ x \end{pmatrix}, \qquad
b_{aug}= \begin{pmatrix} b \\ d \\ 0 \end{pmatrix}.
\end{equation}
Applying the analysis of Section~5 of \cite{sctu:2021b} to this problem yields
Algorithm~\ref{alg:augmented LSE 2}. In exact arithmetic, the main difference between the work required by
Algorithms~\ref{alg:updating_LSE} and \ref{alg:augmented LSE 2}
is that the former involves an additional solve with $RP^T$. For both algorithms, $K$ is
independent of $b$ and $d$.
\begin{algorithm}[htbp]
\caption{\label{alg:augmented LSE 2}
Solve the LS problem (\ref{eq:ls})--(\ref{eq:constraints})
with $C$ having full row rank using the 3-block augmented system (\ref{eq:LSE aug})}
\medskip
\begin{algorithmic}[1]
\State Compute the sparse QR factorization
$\begin{pmatrix} AP & b  \end{pmatrix} = Q\begin{pmatrix} R &f \\ 0 & g \end{pmatrix}$.
\State Solve $P R^TK^T = C^T$  for $K^T \in \Re^{n \times p}$.
\State Compute the minimum-norm solution of $Ku = d - Kf$.
\State Solve $R P^T x = f + u$ for $x$.
\end{algorithmic}
\end{algorithm}

\subsection{Augmented regularized normal equations}

The next approach weights the constraints and uses a regularization parameter within an augmented system formulation
and then aims to balance these two modifications. Consider the weighted least squares problem (WLS)
\begin{equation}\label{eq:weighting}
\min_{x}\left\| A_{\gamma} x_{\gamma} - b_{\gamma} \right\|^2
\;\; \mbox{with} \;\;
A_{\gamma} =  \begin{pmatrix} A \\ \gamma  C  \end{pmatrix}, \;\;
b_{\gamma} = \begin{pmatrix} b \\ \gamma  d \end{pmatrix},
\end{equation}
for some very large $\gamma$ ($\gamma \gg 1$). Let
$x_{LSE}$ be the solution of the LSE problem (\ref{eq:ls})--(\ref{eq:constraints}). Then because
\[ \lim_{\gamma \rightarrow \infty} x_{\gamma} = x_{LSE}, \]
the WLS problem can be used to approximately solve the LSE problem \cite{laha:74}.
An obvious solution method is to solve the normal equations for (\ref{eq:weighting})
\[
H_\gamma x = A_\gamma^T A_\gamma x = (A^TA + \gamma^2 C^TC) x = A^T b + \gamma^2 C^T d= A_\gamma^T b_\gamma.
\]
The appeal of this is that no special methods are required: software for solving
standard normal equations can be used. However, for very large values
of the parameter $\gamma$, the normal matrix $H_\gamma$
becomes extremely ill-conditioned; this is discussed
in Section 4 of \cite{bjdu:80}, where it is shown that the method of
normal equations can break down if $\gamma > \epsilon^{-1/2}$ ($\epsilon$ is
the machine precision). Furthermore, if $C$ contains dense rows then $H_\gamma$ will be dense.

Another possibility is to use the regularized normal equations
\begin{equation}\label{eq:WLS normal reg}
(H_\gamma + \omega^2 I) x = A_\gamma^T b_\gamma,
\end{equation}
where $\omega > 0$ is a regularization parameter \cite{zhgo:2015}. Solving (\ref{eq:WLS normal reg})
is equivalent to solving the $(m + p +n) \times (m + p +n)$ augmented regularized normal equations
\begin{equation}\label{eq:WLS normal aug}
{\mathcal A}(\omega,\gamma)
\begin{pmatrix} y \\ x  \end{pmatrix} =
\begin{pmatrix} b_\gamma \\ 0  \end{pmatrix}, \qquad
{\mathcal A}(\omega,\gamma) =
\begin{pmatrix} \omega I & A_\gamma \\ A_\gamma^T & -\omega I \end{pmatrix},
\end{equation}
where $y = \omega^{-1}(b_\gamma - A_\gamma x)\in \Re^{m+p}$.
The spectral condition number of (\ref{eq:WLS normal aug}) is
\[
 \textrm{cond}({\mathcal A}(\omega,\gamma)) = \sqrt{\textrm{cond}(H_\gamma + \omega^2 I)}
\]
and Saunders~\cite{saun:95a} shows that $\textrm{cond}({\mathcal A}(\omega,\gamma))
\approx \|A_{\gamma}\|/\omega$ regardless of the condition of $A_{\gamma}$.
Thus using (\ref{eq:WLS normal aug}) potentially gives
a significantly more accurate approximation to the pseudo solution $x = A_\gamma^+ b_\gamma$ (where $(.)^+$ denotes
the Moore-Penrose pseudo inverse of a matrix) compared to the
approximation provided by solving (\ref{eq:WLS normal reg}).
In \cite{zhda:2012}, the parameters are set to $\omega = 10^{-q}$ and $\gamma = 10^q$, where
$$q = \min \{ k: 10^{-2k} \le \nu^{-t} \}.$$
Here $t$-bit floating-point arithmetic with base $\nu$ is used.

Rewriting (\ref{eq:WLS normal aug}) using (\ref{eq:weighting})
and a conformal partitioning of $y$ gives
\begin{equation}\label{eq:WLS 3 block}
\begin{pmatrix} \omega I & 0 & A \\  0 & \omega I & \gamma C \\  A^T & \gamma C^T  & -\omega I \end{pmatrix}
 \begin{pmatrix} y_s \\ y_c \\ x  \end{pmatrix} =
\begin{pmatrix} b \\ \gamma  d \\  0  \end{pmatrix}.
\end{equation}
This system can be solved as in \cite{sctu:2020a} using a modified version of
Algorithm~\ref{alg:augmented LSE 2}. Or, eliminating $y_s$ and setting $\omega \gamma = 1$, yields
\begin{equation}\label{eq:WLS 2 block}
\begin{pmatrix} -H(\omega) & C^T \\  C & \omega^2 I \end{pmatrix}
 \begin{pmatrix} x \\ y_c   \end{pmatrix} =
\begin{pmatrix} -A^Tb \\ d  \end{pmatrix},
\qquad H(\omega) = A^TA + \omega^2 I.
\end{equation}
We can solve this system using a QR factorization of $\begin{pmatrix}A \\ \omega I \end{pmatrix}$ and modifying
Algorithm~\ref{alg:updating_LSE}. Or, ignoring the block structure,
we can treat it as a sparse symmetric indefinite linear system
and compute an $LDL^T$ factorization
(with $L$ unit lower triangular and $D$ block diagonal with blocks of size 1 and 2)
using a sparse direct solver such as {\tt HSL\_MA97} \cite{hosc:2013c} that incorporates
pivoting for stability with a  sparsity-preserving ordering. This factorization
would have to be recomputed for each new set of constraints. Alternatively,
a block signed Cholesky factorization of (\ref{eq:WLS 2 block}) can be used, that is,
\begin{equation*}
\begin{pmatrix} -H(\omega) & C^T \\  C & \omega^2 I \end{pmatrix} =
\begin{pmatrix} L &  \\  B & \;L_{\omega} \end{pmatrix}
\begin{pmatrix} -I &  \\  & \;I \end{pmatrix}
\begin{pmatrix} L^T & B^T \\   & L_{\omega}^T \end{pmatrix},
\end{equation*}
where
\[ H(\omega) = LL^T, \quad LB^T = -C^T \quad \mbox{and}
\quad S = \omega^2 I + BB^T = L_{\omega}L_{\omega}^T. \]
We then obtain Algorithm~\ref{alg:weighted augmented LSE}.
Note that $B$ need not be computed explicitly. Rather,
the Schur complement $S$ may be computed using
$\omega^2 I + CL^{-T}L^{-1}C^T$, and $w=Bz$
may be computed by solving $Lv = z$ and then setting $w=-Cv$, and $w= -B^Ty_c$ may be obtained
by solving $Lw = C^Ty_c$.
\begin{algorithm}[htpb]
\caption{\label{alg:weighted augmented LSE}
Given $\omega > 0$, solve the augmented system (\ref{eq:WLS 2 block}) using Cholesky factorizations.}
\medskip
\begin{algorithmic}[1]
\State Compute the sparse Cholesky factorization $H(\omega) = LL^T$.
\State Solve $L z = A^Tb$.
\State Solve $LB^T = -C^T$.
\State Form the symmetric positive definite Schur complement $S = \omega^2 I + BB^T$ and factorize it $S = L_{\omega}L_{\omega}^T$.
\State Solve $L_{\omega} v = d + Bz$ then solve $L^T_{\omega} y = v$.
\State Solve $L^T x = z - B^T y_c$.
\end{algorithmic}
\end{algorithm}

\begin{table}[htbp]
\caption{Results for the augmented regularized normal equations approach  (Algorithm~\ref{alg:weighted augmented LSE}) for
problems sctap1-2r,  south31, and deter3\_20 using a range of values of $\omega$.
$iters$ is the number of GMRES iterations. The computed $\|x\|$ and $\|r\|$
are consistent for both approaches.
\label{T:res weighted aug}
}
\begin{center}
\begin{tabular}{l r  rrl cl}   \hline
\multicolumn{1}{c} { } & & & &
\multicolumn{1}{c} {Algorithm~\ref{alg:weighted augmented LSE}} &
\multicolumn{2}{c} {Preconditioned GMRES} \\
\multicolumn{1}{c} {Identifier} &
   \multicolumn{1}{c}{$\omega$} &
    \multicolumn{1}{c}{$\|x\|$} &
    \multicolumn{1}{c}{$\|r\|$} &    \multicolumn{1}{c}{$\|r_c\|$} &
       \multicolumn{1}{c}{$iters$} &
    \multicolumn{1}{c}{$\|r_c\|$} \\
    \hline

sctap1-2r &   $1.0 \times 10^{-2}$ &   1.441$\times 10^2$ &   1.911$\times 10^2$ &   5.073$\times 10^{-1}$
& 6 &   5.073$\times 10^{-1}$ \\

 &   $1.0 \times 10^{-3}$ &   1.646$\times 10^2$ &   2.067$\times 10^2$ &   7.381$\times 10^{-3}$
 & 6 &   7.381$\times 10^{-3}$ \\

 &   $1.0 \times 10^{-4}$ &   1.649$\times 10^2$ &   2.070$\times 10^2$ &   7.419$\times 10^{-5}$
 & 6 &   7.416$\times 10^{-5}$ \\

 &   $1.0 \times 10^{-5}$ &   1.649$\times 10^2$ &   2.070$\times 10^2$ &   7.711$\times 10^{-7}$
 & 6 &   7.417$\times 10^{-7}$ \\

 &   $1.0 \times 10^{-6}$ &   1.649$\times 10^2$ &   2.070$\times 10^2$ &   1.148$\times 10^{-7}$
 &  2  &   7.417$\times 10^{-9}$ \\

 &   $1.0 \times 10^{-7}$ &   1.649$\times 10^2$ &   2.070$\times 10^2$ &   1.081$\times 10^{-7}$
 & 6 &   7.416$\times 10^{-11}$ \\

 &   $1.0 \times 10^{-8}$ &   1.649$\times 10^2$ &   2.070$\times 10^2$ &   1.201$\times 10^{-7}$
 & 6 &   7.642$\times 10^{-13}$ \\

 &   $1.0 \times 10^{-9}$ &   1.649$\times 10^2$ &   2.070$\times 10^2$ &   1.295$\times 10^{-7}$
 & 6 &   4.095$\times 10^{-13}$ \\

 \hline
south31   &   $1.0 \times 10^{-2}$ &   2.749$\times 10^{1}$ &   1.881$\times 10^{2}$ &   8.341$\times 10^{-5}$
&311 &   8.341$\times 10^{-5}$ \\

   &   $1.0 \times 10^{-3}$ &   2.749$\times 10^{1}$ &   1.881$\times 10^{2}$ &   8.338$\times 10^{-7}$
   &337 &   7.338$\times 10^{-7}$ \\

   &   $1.0 \times 10^{-4}$ &   2.749$\times 10^{1}$ &   1.881$\times 10^{2}$ &   8.338$\times 10^{-9}$
   &352&   8.339$\times 10^{-9}$ \\

   &   $1.0 \times 10^{-5}$ &   2.749$\times 10^{1}$ &   1.881$\times 10^{2}$ &   8.312$\times 10^{-11}$
   &354 &   8.847$\times 10^{-11}$ \\

   &   $1.0 \times 10^{-6}$ &   2.749$\times 10^{1}$ &   1.881$\times 10^{2}$ &   1.017$\times 10^{-12}$
   &354 &   1.057$\times 10^{-11}$ \\

   &   $1.0 \times 10^{-7}$ &   2.749$\times 10^{1}$ &   1.881$\times 10^{2}$ &   1.840$\times 10^{-13}$
   &354 &   1.070$\times 10^{-11}$ \\

   &   $1.0 \times 10^{-8}$ &   2.749$\times 10^{1}$ &   1.881$\times 10^{2}$ &   1.294$\times 10^{-13}$
   &354 &   1.073$\times 10^{-11}$ \\

   &   $1.0 \times 10^{-9}$ &   2.749$\times 10^{1}$ &   1.881$\times 10^{2}$ &   6.768$\times 10^{-14}$
   &354 &   1.076$\times 10^{-11}$ \\

  \hline

   deter3\_20
   &   $1.0 \times 10^{-2}$ &   1.218$\times 10^3$ &   1.227$\times 10^2$ &   6.877$\times 10^{-4}$ & 34  &   6.877$\times 10^{-4}$ \\
   &   $1.0 \times 10^{-3}$ &   1.585$\times 10^3$ &   1.220$\times 10^2$ &   6.834$\times 10^{-6}$ & 36  &   6.834$\times 10^{-6}$ \\
   &   $1.0 \times 10^{-4}$ &   1.589$\times 10^3$ &   1.220$\times 10^2$ &   6.834$\times 10^{-8}$ & 36  &   6.834$\times 10^{-8}$ \\
   &   $1.0 \times 10^{-5}$ &   1.589$\times 10^3$ &   1.220$\times 10^2$ &   6.834$\times 10^{-10}$ & 36  &   6.831$\times 10^{-10}$\\
   &   $1.0 \times 10^{-6}$ &   1.589$\times 10^3$ &   1.220$\times 10^2$ &   6.935$\times 10^{-12}$ & 36  &   6.718$\times 10^{-12}$\\
   &   $1.0 \times 10^{-7}$ &   1.589$\times 10^3$ &   1.220$\times 10^2$ &   1.138$\times 10^{-12}$ & 36  &   1.111$\times 10^{-12}$\\
   &   $1.0 \times 10^{-8}$ &   1.589$\times 10^3$ &   1.220$\times 10^2$ &   1.433$\times 10^{-12}$ & 36  &   1.043$\times 10^{-12}$ \\
   &   $1.0 \times 10^{-9}$ &   1.589$\times 10^3$ &   1.220$\times 10^2$ &   1.350$\times 10^{-12}$ & 36  &   1.372$\times 10^{-12}$ \\

   \hline
\end{tabular}
\end{center}
\end{table}
Results for Algorithm~\ref{alg:weighted augmented LSE} for three of our test problems
using  a range of values of $\omega$
are given in Table~\ref{T:res weighted aug}. Note that here $\|r_c\|$
is computed using $r_c = d - C x$ (rather than using $r_c = \omega*y_c$).
We see that, provided $\omega$ is sufficiently small,  the values
of  $\|x\|$ and  $\|r\|$
are consistent with those given in Table~\ref{T:test problems}.

By replacing the Cholesky factorization of $H(\omega)$ by an incomplete factorization $H(\omega) \approx \tilde{L}\tilde{L}^T$,
we can obtain a preconditioner for solving (\ref{eq:WLS 2 block}).
In particular, the right-preconditioned system is
\begin{equation} \label{eq:aug_pre}
\begin{pmatrix} -H(\omega) & C^T \\  C & \omega^2 I \end{pmatrix}
       M^{-1}\begin{pmatrix} w \\ w_c  \end{pmatrix} = \begin{pmatrix} -A^T b \\ d \end{pmatrix},
 \qquad M\begin{pmatrix} x \\ y_c  \end{pmatrix} = \begin{pmatrix} w \\ w_c  \end{pmatrix},
   \end{equation}
   and we can take the preconditioner in factored form to be
\begin{equation} \label{eq:aug_precon}
M
=
\begin{pmatrix} \tilde{L} \\[0.1cm]    \tilde{B} & \;\;\;I \end{pmatrix}
\begin{pmatrix} -I & \\[0.1cm]  &  \tilde{S}_d \end{pmatrix}
\begin{pmatrix} \tilde{L}^T & \tilde{B}^T \\[0.1cm]    & I \end{pmatrix},
\end{equation}
with
\[
\tilde{L}\tilde{B}^T = -C^T \quad \mbox{and} \quad
\tilde{S} = \omega^2 I + \tilde{B} \tilde{B}^T.
\]
As the preconditioner (\ref{eq:aug_precon}) is indefinite, it needs to be
used with a general nonsymmetric iterative
method such as GMRES~\cite{sasc:86}. A positive definite preconditioner
for use with MINRES~\cite{pasa:75} can be obtained by replacing $-I$ in (\ref{eq:aug_precon}) by $I$. MINRES has the important advantage of only requiring three
vectors of length equal to the size of the linear system. 
GMRES results are included in Table~\ref{T:res weighted aug}. 
The GMRES convergence tolerance is
taken to be $10^{-11}$. We see that the GMRES iteration count is essentially
independent of $\omega$. We also ran MINRES with the same
settings and, while the iteration
counts were again insensitive to $\omega$, they were significantly 
greater than for GMRES. For problems sctap1-2r, south31 and deter3\_20 the counts
were 17, 772 and 56, respectively  ($\omega = 1.0\times 10^{-5}$).

\begin{table}[htbp]
\caption{Convergence results for problems sctap1-2r with $\omega=1.0\times 10^{-8}$
and stormg2-8\_20 with $\omega=1.0\times 10^{-6}$.
$tol$ and $iters$ are the convergence tolerance and the iteration count for GMRES.
\label{T:res stormg2-8}
}
\begin{center}
\begin{tabular}{l cl cl}   \hline
&
   \multicolumn{2}{c}{sctap1-2r} &
      \multicolumn{2}{c}{stormg2-8\_20} \\

      \multicolumn{1}{c}{$tol$} &
    \multicolumn{1}{c}{$iters$} &
    \multicolumn{1}{c}{$\|r_c\|$ } &
    \multicolumn{1}{c}{$iters$} &
    \multicolumn{1}{c}{$\|r_c\|$ } \\
    \hline
   $1.0 \times 10^{-6}$ &   2   &   1.669$\times 10^{-6}$  & 130 &   1.423$\times 10^{-7}$ \\
   $1.0 \times 10^{-7}$ &   3   &   6.046$\times 10^{-8}$  & 134 &   5.364$\times 10^{-9}$ \\
   $1.0 \times 10^{-8}$ &   3   &   6.046$\times 10^{-8}$  & 141 &   1.434$\times 10^{-9}$ \\
   $1.0 \times 10^{-9}$ &   4   &   1.897$\times 10^{-8}$  & 146 &   1.246$\times 10^{-10}$ \\
   $1.0 \times 10^{-10}$ &   4  &   1.897$\times 10^{-8}$  & 149 &   9.067$\times 10^{-11}$ \\
   $1.0 \times 10^{-11}$ &   6  &   7.642$\times 10^{-13}$ & 156 &   1.314$\times 10^{-10}$ \\
   $1.0 \times 10^{-12}$ &   6  &   7.642$\times 10^{-13}$ & 161 &   5.220$\times 10^{-11}$ \\
   $1.0 \times 10^{-13}$ &   6  &   7.642$\times 10^{-13}$ & 190 &   4.972$\times 10^{-12}$ \\
   $1.0 \times 10^{-14}$ &   7  &   9.136$\times 10^{-13}$ & 217 &   4.974$\times 10^{-12}$ \\
  \hline
\end{tabular}
\end{center}
\end{table}

%
%


Our findings in Section~\ref {sec:aug system} suggest that,
if we require the constraints to be solved with a small residual,
then an augmented system based approach combined with a QR factorization performs better (in terms of $\|r_c\|$) than
combining it with regularization and a Cholesky factorization.
Unfortunately, QR factorizations are more expensive and while strategies for computing
incomplete orthogonal factorizations for use in building preconditioners have been proposed
(see, for instance, \cite{bdw:01,bady:09,bayi:09,jeaj:84,lisa:06,pdw:05,wagb:97}),
the only available software is the MIQR package of Li and Saad~\cite{lisa:06}
(probably because developing high quality implementations is non-trivial).
In their study of preconditioners for LS problems, Gould and Scott~\cite{gosc:2015b,gosc:2017}
found that MIQR generally performed less well than incomplete Cholesky factorization preconditioners
and so is not considered here.

We have made the implicit assumption that $A$ is sparse. However, it is straightforward to extend the
augmented system-based approaches to the more general case that $A$ contains
rows that are dense. For example, if $A$ is permuted and partitioned as
\[  A = \begin{pmatrix} A_1 \\  A_2   \end{pmatrix},\]
where $A_1$ is sparse and $A_2$ is dense, then using a conformal partitioning of $y_s$ and of $b$, (\ref{eq:WLS 2 block}) can be replaced by
the augmented system
\[
\begin{pmatrix}  -H_{1}(\omega) & C_d^T \\[0.15cm]
                  C_d & \omega^2 I \end{pmatrix}
 \begin{pmatrix} x \\[0.15cm] y_d   \end{pmatrix} =
\begin{pmatrix} -A_1^Tb_1 \\[0.15cm] d  \end{pmatrix}
\]
with
\[
\quad H_{1}(\omega) = A_1^TA_1 + \omega^2 I,
\quad  y_d = \begin{pmatrix} y_c \\ y_2   \end{pmatrix},
\quad C_d = \begin{pmatrix} C \\ \omega A_2   \end{pmatrix},
 \quad d = \begin{pmatrix} d \\ \omega b_2   \end{pmatrix}.
\]

Finally, we remark that, if we use the 3-block form (\ref{eq:WLS 3 block})
then we can follow \cite{sctu:2020a}, which in turn generalises
the work of Carson, Higham and Pranesh~\cite{cahp:2020a}, and obtain an
augmented system approach with multi-precision refinement. This has the potential
to reduce the computational cost in terms of time and/or memory, thus allowing
larger problems to be solved.

\section{Conclusions} \label{sec:conclusions}
We have considered a number of approaches for solving large-scale LSE problems
in which the constraints may be dense.
Our main findings can be summarized as follows:
\begin{itemize}
 \item The classical null-space method relies on computing a null-space basis matrix $Z$ for the ``wide''
 constraint matrix $C$ such that  $Z^T A^T A Z$ is sparse. In recent work \cite{sctu:2021b},
 we proposed how this can be achieved using a method based on a QR factorization
 of $C$ with threshold
 pivoting. This is not straightforward to implement.
 Furthermore, our numerical experiments show that, in some cases, the norm $\|r_c\|$ of the
 constraints residual can be larger than for other approaches considered in this study. Thus,
 although in some contexts null-space approaches are popular, we do not recommend the
 strategy of \cite{sctu:2021b} for LSE problems.

 \item The direct elimination approach
 couples the constraint matrix and the LS matrix, leading to a sparse-dense transformed
 least squares problem. Existing direct or iterative methods can be used to solve the transformed problem
 and our experiments found the computed constraint residuals are small.
 The approach can be used for problems for which $A$ (as well as $C$)
contains a small number of dense rows. A weakness is that, if solving a sequence of problems in which either
$A$ or $C$ is fixed, the coupling of the two blocks in the solution process means that it must be restarted.
Furthermore, the number of dense rows in the transformed problem can be relatively large, making it expensive to solve.

 \item There are several options for using an augmented system formulation.
 This can be solved using standard building blocks, such as a sparse QR factorization, a sparse symmetric indefinite linear solver,
 or a block sparse Cholesky factorization.
 An attraction of each of these is that existing
 ``black box'' solvers can be exploited, thereby greatly reducing the effort required in developing
 robust and efficient implementations.  The augmented system formulation can be generalised to handle dense rows in $A$
 and offers the potential for mixed-precision computation.
 Moreover, an incomplete Cholesky factorization can be used as a
 preconditioner with a Krylov subspace solver.

 \item In the case of a series of LSE problems in which only the constraints
 change, both the null-space and direct elimination approaches have the disadvantage that the computation
 must be redone for each new set of constraints. For the augmented system approaches,
 a significant amount of work can be reused from the first problem in the sequence
 when solving subsequent problems.

\end{itemize}
\smallskip
Finally, we observe that there is a lack of iterative methods and preconditioners that can be used to extend the
 size of LSE problems that can be solved. We have shown that using an incomplete factorization
 within a block factorization of an augmented system can be effective, but
 most current incomplete factorizations
 that result in efficient preconditioners are serial in nature and not able
 to tackle extremely large problems (but see \cite{ancd:2018,hsth:2018} for novel
 approaches that are designed to exploit parallelise).
 Addressing  the lack of iterative approaches is a challenging subject for future work.

\begin{acknowledgements}
We are grateful to Professor Michael Saunders and an anonymous reviewer
for their constructive comments that have led to many improvements in the presentation
of this paper.
\end{acknowledgements}

\bibliographystyle{plain}      

\bibliography{btbook}   

\end{document}